\documentstyle[12pt]{article}

\def\fr{\frac}
\def\la{\lambda}

\def\Si{\Sigma}

\def\de{\delta}

\def\a{a^{\dagger}}

\newcommand{\bd}{\begin{displaymath}}
\newcommand{\ed}{\end{displaymath}}
\newcommand{\bb}{\begin{equation}}
\newcommand{\ee}{\end{equation}}
\newcommand{\bea}{\begin{eqnarray}}
\newcommand{\eea}{\end{eqnarray}}

\begin{document}
\baselineskip 1.5 \baselineskip

%% Title Page 

\vspace{.2cm}

\begin{center}
\Large {\bf Coherent States of $gl_q(2)$-covariant Oscillators }
\\[1cm]
\large W.-S.Chung \\[.3cm]
\normalsize  
Theory Group, Department of Physics and Research Institute of Natural
Science,  \\
\normalsize  Gyeongsang National University,   \\
\normalsize   Jinju, 660-701, Korea
\end{center}

\vspace{0.5cm}
\begin{abstract}
In this paper two types of coherent states of $gl_q(2)$-covariant
oscillators are
investigated.
\end{abstract}

\setcounter{page}{1}
\section{}

Quantum                                                               
groups                                                                      
or                        
q-deformed                    
Lie         
algebra                
implies   
some    
specific 
deformations 
of                    classical                      Lie           algebras.

From           
a             
mathematical      
point   
of  
view,     
it 
is  
a 
non-commutative                                                  
associative                                                                 
Hopf         
algebra.         
The       
structure        
and 
representation 
theory    of                                   quantum               groups  
have                 
been             
developed        
extensively  
by  
Jimbo    
[1] 
and 
Drinfeld                                                                
[2].                                                                        
            
The                                   
q-deformation                    
of        
Heisenberg       
algebra   
was     
made  
by 
Arik and Coon [3], Macfarlane [4] and Biedenharn [5].
Recently                                                               
there                                               
has                           
been                    
some              
interest               
in    
more     
general   
deformations 
involving                                                                 
an                                                                          
arbitrary                      
real                           
functions     
of           
weight    
generators  
and   
including 
q-deformed algebras as a special case [6-10].

Recently V.Spiridonov [11] found the new coherent states of the q-Weyl algebra
$a\a -q \a a =1, 0<q<1$ which are defined as eigenstates of the operator $\a$ .

The purpose of this paper is to find 
such coherent states for $gl_q(2)$-covariant oscillator algebra.

\def\a{a^{\dagger}}
\def\qq{\fr{1}{\sqrt{q}}}
\def\q{\sqrt{q}}

The $gl_q(2)$-covariant oscillator algebra is defined as
\begin{displaymath}
 \a_1 \a_2=\q\a_2\a_1
\end{displaymath}
\begin{displaymath}
a_1a_2=\qq a_2a_1
\end{displaymath}
\begin{displaymath}
a_1\a_2=\q\a_2a_1
\end{displaymath}
\begin{displaymath}
a_2\a_1=\q\a_1a_2
\end{displaymath}
\begin{displaymath}
a_1\a_1=1+q\a_1a_1+(q-1) \a_2a_2
\end{displaymath}
\bb
a_2\a_2=1+q\a_2a_2
\ee
Throughtout, $()^{\dagger}$ denotes the hermitian conjugate of $()$.
By the $gl_q(2)$-covariance of the system, it is meant that the linear
trnsformations
\bd
\matrix{ a & b \cr c & d \cr}\matrix{ a_1 \cr a_2\cr }=\matrix{ a_1^{\prime}\cr
a_2^{\prime}\cr}
\ed
\bb
\matrix{ \a_1 & \a_2} \matrix{a &b \cr c& d \cr}=\matrix{(\a_1)^{\prime}&
(\a_2)^{\prime}}
\ee
leads to the same commutation relations (1) for 
$( a_1^{\prime},(\a_1)^{\prime})$  and
$(a_2^{\prime},(\a_2)^{\prime})$ when the matrix $ \matrix{a &b\cr c & d
\cr}$ belongs to
$gl_q(2)$.
The relation between the entries of this matrix obeys the following
commutation relation
\bd
ab=\q ba~~~~cd=\q dc
\ed
\bd
ac=\q ca~~~~bd=\q db
\ed
\bb
bc=cb ~~~~ad-da=(\q-\qq)bc
\ee
It should be noted that the particular coupling between the two modes 
 is completely 
dictated by the required $gl_q(2)$-covariance.

The Fock space representation of the algebra (1) can be easily constructed by
introducing the hermitian number operators $\{ N_1, N_2\}$ obeying
\bb
[N_i, a_j]=-\de_{ij}a_j~~~[N_i,\a_j]=\de_{ij}\a_j,~~~(i,j=1,2)
\ee
Let $|0,0>$ be the unique ground state of this system satisfying
\bb
N_i|0,0>=0,~~~a_i|0,0>=0,~~~(i,j=1,2)
\ee
and $\{ |n,m>|n,m=0,1,2,\cdots\}$ be the set of the orthogonal number
eigenstates
\bd
N_1|n,m>=n|n,m>,~~~~N_2|n,m>=m|n,m>
\ed
\bb
<n,m|n^{\prime},m^{\prime}>=\de_{n n^{\prime}}\de_{m m^{\prime}}
\ee
From the algebra (1) the representation is given by

\def\n{|n,m>}
\def\nb{|n,m-1>}
\def\na{|n-1,m>}
\begin{displaymath}
a_1\n =\sqrt{q^{m}[n]}\na,~~~
a_2\n=\sqrt{[m]}\nb
\end{displaymath}
\bb
\a_1\n =\sqrt{q^{m}[n+1]}|n+1,m>,~~~
\a_2\n=\sqrt{[m+1]}|n,m+1>
\ee
where the q-number $[x]$ is defined as
\bd
[x]=\fr{q^x-1}{q-1}
\ed
The general eigenstates $\n$ is obtained  by applying $\a_2$ m times after
applying $\a_1$ 
n times.
\bb
\n =\fr{(\a_2)^m (\a_1)^n }{\sqrt{[n]![m]!}}|0,0>
\ee
where 
\bd
[n]!=[n][n-1]\cdots[2][1],~~~[0]!=1
\ed
The coherent states for algebra (1)  are usually defined as 
\bd
a_1|z_1,z_2>_-=z_1|z_1,\q z_2>_-
\ed
\bb
a_2|z_1,z_2>_-=z_2|z_1,z_2>_-
\ee
It can be easily checked that the coherent state satisfies $a_1a_2 =\qq a_2
a_1$.
Solving the relation (9) we have
\bb
|z_1,z_2>_-=c(z_1,z_2)\Si_{n,m=0}^{\infty}\fr{z_1^nz_2^m}{\sqrt{[n]![m]!}}\n
\ee
Using the eq.(8) we can rewrite eq.(9) as
\bb
|z_1,z_2>_-=c(z_1,z_2)e_q(z_1\a_1)e_q(z_2 \a_2)|0,0>
\ee
In order to obtain the normaized coherent states, we should impose the condition
$~{}_<z_1,z_2|z_1.z_2>_-=1$. Then the normalized coherent states are given by
\bb
|z_1,z_2>_-=\fr{1}{\sqrt{e_q(|z_1|^2)e_q(|z_2|^2)}}e_q(z_1\a_1)e_q(z_2
\a_2)|0,0>
\ee
where
\bd
e_q(x)=\Si_{n=0}^{\infty}\fr{x^n}{[n]!}
\ed
is a q-deformed exponential function.

\section{}
The purpose of this section is to obtain another type of 
coherent states for algebra (1).
In order to do so , it is convenient to introduce
the two commuting operators as follows
\begin{displaymath}
H=\a_1a_1+\a_2a_2-\nu,~~~T=\a_2a_2-\nu
\end{displaymath}
where
\begin{displaymath}
\nu=\frac{1}{1-q}
\end{displaymath}
and $H$ is a hamiltonian and $T$ is an integral motion for $H$, i.e. an 
independent operator commuting with the hamiltonian $H$.
From the algebra (1) we can easily check that $[H,T]=0$.
Then the commutation relation between two commuting operators and mode
operators are given 
by
\begin{displaymath}
H\a_1=q\a_1H,~~~T\a_2=q\a_2T
\end{displaymath}
\begin{displaymath}
H\a_2=q\a_2H,~~~~T\a_1=\a_1T
\end{displaymath}
Acting the two commuting operators on the eigenstates gives
\bb
H\n =E_{nm}\n~~~T\n =t_m\n
\ee
where
\bb
E_{nm} =-\fr{q^{m+n}}{1-q},~~~~t_m=-\fr{q^m}{1-q}
\ee
It is worth noting that the energy spectrum is degenerate and negative for
$0<q<1$.
The degenerated states are splited by the eigenvalues of operator $T$, but
the physical interpretation of operator $T$ is not clear to the author.

As was  noticed in  ref[11], for  the positive  energy states  it is  not
$a_1(  a_2)$ but 
$\a_1(\a_2)$ that
play a role of the lowering operator:

\begin{displaymath}
H|\lambda q^{n},\mu q^{m}>,~~n,m=0,\pm 1 ,\pm 2 , \cdots
=\lambda q^{n+m}|\lambda q^{n},\mu q^{m}>
\end{displaymath}
\begin{displaymath}
T|\lambda q^{n},\mu q^{m}>
=\mu q^{m}|\lambda q^{n},\mu q^{m}>
\end{displaymath}
\begin{displaymath}
\a_1|\lambda q^{n},\mu q^{m}>
=\sqrt{\lambda q^{n+1}-\mu q^{m}}|\lambda q^{n+1},\mu q^{m}>
\end{displaymath}
\begin{displaymath}
\a_2|\lambda q^{n},\mu q^{m}>
=\sqrt{\mu q^{m+1}+\nu}|\lambda q^{n+1},\mu q^{m+1}>
\end{displaymath}
\begin{displaymath}
a_1|\lambda q^{n},\mu q^{m}>
=\sqrt{\lambda q^{n}-\mu q^{m}}|\lambda q^{n-1},\mu q^{m}>
\end{displaymath}
\bb
a_2|\lambda q^{n},\mu q^{m}>
=\sqrt{\mu q^{m}+\nu}|\lambda q^{n-1},\mu q^{m-1}>
\ee
where $ \lambda, \mu >0$ are arbitrary chosen eigenvalues of $H$ and $T$.
This representation is a nonhighest weight representation of the
$gl_q(2)$-covariant
oscillator algebra. And this representation has no classical analogue because
it is not defined for $q \rightarrow 1$.
Due to this fact, it is natural to define coherent states 
corresponding to the representation (15)
 as the eigenstates of $\a_1$ and $\a_2$:

\def\s{\Sigma}
\def\fr{\frac}
\def\l{\lambda}
\def\m{\mu}
\def\z{|z_1,z_2>_+}
\def\zz{|z_1,\qq z_2>_+}
\bb
\a_1\z =z_1 \zz,~~~\a_2\z=z_2\z
\ee
Because  the representation  (15)  depend on  two free  paprameter  $ \la$
and  $\m$, the 
coherent states
$\z$ can take different forms.

If we assume that the positive energy states are normalizable, 
i.e.$~~~$    $<\la    q^n,    \m    q^m|\la    q^{n^{\prime}},\m
q^{m^{\prime}}>=\de_{n 
n^{\prime}}\de_{mm^{\prime}}$, 
and form exactly one series for some 
fixed $\la$ and $\m$, then we can write
\bb
|z_1,z_2>_+=\s_{n,m=-\infty}^{\infty} c_{nm}(z_1,z_2)|\lambda q^{n},\mu q^{m}>
\ee
Inserting eq.(17) into eq.(16), we find

\bea
|z_1,z_2>_+
&=&C[ \s_{n=1}^{\infty}\s_{m=1}^{\infty}
c_{nm}^1z_1^{m-n}z_2^{-m}|\la q^n, \m q^m>\cr
&+&\s_{n=1}^{\infty}\s_{m=0}^{\infty}
c_{nm}^2z_1^{-m-n}z_2^{m}|\la q^n , \m q^{-m}>\cr
&+&\s_{n=0}^{\infty}\s_{m=1}^{\infty}
c_{nm}^3z_1^{m+n}z_2^{-m}|\la q^{-n},\m q^m>\cr
&+&\s_{n=0}^{\infty}\s_{m=0}^{\infty}
c_{nm}^4z_1^{n-m}z_2^{m} |\la q^{-n},\m q^{-m}> ]\cr
\eea
where
\bea
c_{nm}^1&=&q^{\frac{ m(m-1)}{4}}
\frac{(\sqrt{-\mu})^n (\sqrt{\nu})^m}{(\sqrt{\lambda})^m}
\sqrt{\frac{(\frac{\lambda}{\mu}q^{1-m};q)_n
(-\frac{\mu}{\nu}q;q)_m}
{(\fr{\m}{\l};q)_m}}\cr
c_{nm}^2&=&q^{\frac{m(m-1)}{4}}\frac{(\sqrt{-\mu})^n
(\sqrt{-\mu})^m}{(\sqrt{\m})^m}
\sqrt{\frac{(\frac{\lambda}{\mu}q^{m+1};q)_n
(\frac{\l}{\mu}q;q)_m}
{(-\fr{\nu}{\m};q)_m}}\cr
c_{nm}^3&=&q^{\fr{nm}{2}+\frac{m(m-1)}{4}+\fr{n(n-1)}{4}}
\frac{(\sqrt{\nu})^m }{(\sqrt{\l})^{n}(\sqrt{\lambda})^m}
\sqrt{\frac{(-\frac{\m}{\nu}q;q)_m}
{(\frac{\mu}{\l}q^{m};q)_{n}
(\fr{\m}{\l};q)_m}}\cr
c_{nm}^4&=&q^{-\fr{nm}{2}+\frac{m(m-1)}{4}+\fr{n(n-1)}{4}}
\frac{(\sqrt{-\mu})^m }{(\sqrt{\l})^{n}(\sqrt{\mu})^{m}}
\sqrt{\frac{(\frac{\l}{\mu}q;q)_m}
{(\frac{\mu}{\l}q^{-m};q)_{n}
(-\fr{\nu}{\m};q)_{m}}}\cr
\eea

If we demand that ${}_+<z_1,z_2|z_1,z_2>_+=1$, we have
\bb
|C|^{-2} =\s_{n=-\infty}^{\infty}\s_{m=-\infty}^{\infty}
\fr{(-)^m q^{-nm +\fr{n(n-1)}{2}+\fr{m(m-1)}{2}}(\fr{\l}{\m}q;q)_m}
{\l^n (\fr{\m}{\l}q^{-m};q)_n(-\fr{\nu}{\m};q)_m}|z_1|^{2(n-m)}|z_2|^{2m}
\ee
If we substitute $ n-m\rightarrow l $ in eq.(20) and use the identity
\bb
(aq^{-m};q)_n=(-a)^mq^{-\fr{m(m+1)}{2}}(q/a;q)_m(a;q)_{n-m}
\ee
we can express the normalization constant $C$ 
in terms of the bilateral q-hypergeometric series [12]:
\bb
|C|^{-2} =
{}_0\psi_1(-\fr{\nu}{\m};q,-\fr{|z_2|^2}{\m})
~{}_0\psi_1(\fr{\m}{\l};q,-\fr{|z_1|^2}{\l})
\ee
where general bilateral q-hypergeometric series is defined by [12]
\bb
{}_r \psi_s \left(\matrix{ a_1,&\cdots,& a_r \cr b_1,& \cdots,& b_s\cr}
;q ,z \right)
=\s_{n=-\infty}^{\infty}
\fr{(a_1;q)_n\cdots (a_r;q)_n}{(b_1;q)_n\cdots (b_s;q)_n}((-)^n
q^{n(n-1)/2})^{s-r}z^n
\ee

We can introduce the nonunitary displacement operator $D(z_1,z_2)$
\bb
\z=D(z_1,z_2)|\l,\m>
\ee
Then the operator $D(z_1,z_2)$ is given by
\bb
D(z_1,z_2)=\fr{{}_0\psi_1(-\nu/\m;q,-z_2a_2/\m){}_0\psi_1(\m/\l ;q,-z_1a_1/\l)}
{\sqrt{
{}_0\psi(-\fr{\nu}{\m};q,-\fr{|z_2|^2}{\m})
{}_0\psi_1(\fr{\m}{\l};q,-\fr{|z_1|^2}{\l})}}
\ee

Formally we can write also
\bd
D(z_1,z_2)=\s_{n,m=-\infty}^{\infty} \left( \fr{\a_2}{z_2}\right)^m \left(
\fr{\a_1}{z_1}\right)^{n-m}
\ed
or
\bea
D(z_1,z_2)&=&\s_{n=1}^{\infty}\s_{m=1}^{\infty} \left(
\fr{\a_2}{z_2}\right)^m \left(
\fr{\a_1}{z_1}\right)^{n-m}                   \cr
&+&\s_{n=1}^{\infty}\s_{m=0}^{\infty}\fr{q^{\fr{m(m-1)}{2}}}{(-\fr{\nu}{\m};
q)_m}   \left( 
\fr{z_2a_2}{\m}\right)^m \left(
\fr{\a_1}{z_1}\right)^{n+m}                   \cr
&+&\s_{n=0}^{\infty}\s_{m=1}^{\infty}\fr{q^{\fr{1}{2}(n+m)(n+m-1)}}{(\fr{\nu
}{\l} q^m;q)_n 
(\fr{\m}{\l};q)_m} \left( \fr{\a_2}{z_2}\right)^m \left(
\fr{z_1a_1}{\l}\right)^{n+m}                   \cr
&+&\s_{n=0}^{\infty}\s_{m=0}^{\infty}\fr{q^{-nm+\fr{n(n-1)}{2}+\fr{1}{2}m(m-
1)}}{(\fr{\mu}
{\l} q^{-m};q)_n (-\fr{\nu}{\m};q)_m} \left( \fr{z_2a_2}{\l}\right)^m \left(
\fr{z_1a_1}{\l}\right)^{n-m}                   \cr
\eea

For some realizations  of $gl_q(2)$-covariant oscillator algebra,  the
states (24) belongs 
to a continuous spectrum.
Thus it is appropriate to consider integrals over both $\l$ and $\m$ for the
expansion of
$\z$ in the basis $|\l,\m>$ instead of sum :
\bb
\z =\int_0^{\infty}\int_0^{\infty}
d\l d\m C(\l,\m,z_1,z_2)|\l,\m>
\ee
Inserting this expression into eq.(16) gives
\bb
\z =\int_0^{\infty}\int_0^{\infty} d\l d\m \fr{\l^d \m^{\epsilon}h(\l)|\l,\m>}
{\sqrt{(q\l/\m;q)_{\infty}(-q\m/\nu;q)_{\infty}}}
\ee
where
\bd
d=-\fr{\ln  q z_1}{\ln q},~~~\epsilon=\fr{\ln \fr{z_1}{q z_2}\sqrt{\nu}}{\ln q}
\ed
where $h(\l)$ is an arbitrary function satisfying $h(q\l)=h(\l)$.

If we impose the normalization condition
\bd
<\l,\m|\sigma,\eta>=\l\m \de(\l-\sigma)\de(\m-\eta)
\ed
then we find that the states are normalizable,

\bb
_+<z_1,z_2|z_1,z_2>_+=1=\int_0^{\infty}\int_0^{\infty}
d\l d\m\fr{\l^{2Re~ d +1}\m^{2Re~ \epsilon +1}|h(\l)|^2}{(q\l/\m;q)_{\infty}
(-q\m/\nu;q)_{\infty}}
\ee
if $h(\l)$ is a bounded function and $ Re~ d >-1, Re~\epsilon >-1$ , or $
|z_1|^2<1, \fr{|z_2|^2}{|z_1|^2} >\nu $.
Expanding $h(\l)$ into the Fourier series, we have an infinite number of
linearly independent ( but not orthogonal ) coherent states of the form
\bb
|z_1,z_2>_{+s} =C(z_1,z_2)
\int_0^{\infty}\int_0^{\infty} 
d\l d\m \fr{\l^{\delta_s} \m^{\epsilon }|\l,\m>}
{\sqrt{(q\l/\m;q)_{\infty}(-q\m/\nu;q)_{\infty}}}
\ee
where
\bd
\delta_s = d + \fr{2 \pi i s }{\ln q }, ~~s=0,\pm 1 ,\cdots
\ed
From  the  normalization   condition  of  the  coherent  states,   we  can
determine  the 
normalization
constatant,
\bb
|C(z_1,z_2)|^{-2} 
=\int_0^{\infty}\int_0^{\infty} 
d\l d\m \fr{\l^{\tau} \m^{\xi}|\l,\m>}
{(q\l/\m;q)_{\infty}(-q\m/\nu;q)_{\infty}}
\ee
where
\bd
\tau =\fr{\ln \fr{1}{q|z_1|^2}}{\ln q},~~
\xi =\fr{\ln \fr{|z_1|^2 \nu}{q^2|z_2|^2}}{\ln q}.
\ed
\section{}
To conculude, in this paper I have constructed two types of coherent states,
one of which is related to the positive energy representation of the
$gl_q(2)$-covariant
oscillator system. These states are shown to be related to the product of two
bilateral basic hypergeometric series. I think that this method will be
applied to
more general case, $gl_q(n)$-covariant oscillator system.
In that case I guess that the new coherent staes associated with the
positive energy 
will be related to the product of $n$ bilateral basic hypergeometric series.
I hope that this problem and its related topics will become clear in the
near future.

\section*{Acknowledgement}
The author want to give many thanks to V.Spiridonov for helpful comment and
advice.
This                   paper                was
supported         by  
the   KOSEF (961-0201-004-2)   
and   the   present   studies    were   supported   by   Basic  
Science 
Research Program, Ministry of Education, 1995 (BSRI-95-2413).

\vfill\eject

\end{document}